\documentclass[12pt,a4paper]{article}
\usepackage{enumerate}
\usepackage{comment}
\usepackage{amsfonts,amsthm,amsmath}
\usepackage[usenames]{color}
\usepackage[pdftex,plainpages=false,pdfpagelabels]{hyperref}

\theoremstyle{plain}
\newtheorem{thm}{Theorem}[section]
\newtheorem{cor}[thm]{Corollary}
\newtheorem{lem}[thm]{Lemma}
\newtheorem{prop}[thm]{Proposition}
\newtheorem{ex}[thm]{Example}
\theoremstyle{remark}
\newtheorem{rem}[thm]{Remark}

\theoremstyle{definition}

\def\RR{\mathbb{R}}

\def\Rm{\mathbb{R}_{\mathrm{max}}}
\def\cX{\mathcal{X}}
\def\cY{\mathcal{Y}}
\def\cD{\mathcal{D}}
\def\cC{\mathcal{C}}
\def\cE{\mathcal{E}}

\def\cT{\mathcal{T}}

\def\Tt{\left(T(t)\right)_{t\ge 0}}
\def\+{\oplus_{\cX}}

\def\cCf{\mathcal{C}_{\mathrm{fin}}}
\def\cDf{\mathcal{D}_{\mathrm{fin}}}

\def\ll1{\mathrm{L}^1}
\def\llp{\mathrm{L}^p}
\def\lli{\mathrm{L}^{\infty}}
\def\cc{\mathrm{C_c}}
\def\cb{\mathrm{C_b}}
\def\cc1{\mathrm{C^1_c}}
\def\c1{\mathrm{C^1}}
\def\supp{\mathrm{supp\,}}
\def\buc{\mathrm{BUC}}

\def\sgn{\mathrm{sgn}}

\title{Semigroups of max-plus linear operators}
\author{Marjeta Kramar Fijav\v{z}, Aljo\v{s}a Peperko and Eszter Sikolya}
\date{\today} 

\begin{document}

\maketitle

\begin{abstract}
We define strongly continuous max-additive and max-plus linear operator semigroups and study their main properties. We present some important examples of such semigroups coming from non-linear evolution equations.
\end{abstract}

\noindent
{\it Math.~Subj.~Classification (2010)}:  47H20, 47J35,15A80. 
\\

\section{Introduction}

The algebraic system of a max-plus algebra and its isomorphic versions (tropical algebra, max algebra) provide 
an attractive way of describing  non-linear problems appearing, e.g., in combinatorial optimisation, mathematical physics, manufacturing and transportation scheduling, information 
technology, discrete event-dynamic systems, DNA analysis, stochastic analysis, etc.~(see \cite{BCOQ92, Bu10,  G92, HOW05, KM, LM05, Mc06} and the references therein). Its usefulness arises from the fact that these non-linear problems become linear when described in the max-plus algebra language. Although the theory was initially developed by algebraists and discrete mathematicians it soon caused interest also in analysis (see  e.g.~\cite{AGN12, KM, LM05, LMS, Mc06}).

Having these facts in mind we introduce the structure of strongly continuous  semigroups of max-plus linear operators. For strongly continuous semigroups of linear operators on Banach spaces, the theory and its applications are already well established and we refer to the classical textbooks \cite{EN00, Go85, Paz83}.  Considerable work has also been done in the non-linear part of the theory, which is much more complicated and branched out then the linear case,  see e.g. ~\cite{B76, BM98,CL71,Mo88,Pa87} and the survey paper \cite{GN14}. Often different approaches are used for different nonlinear  problems arising from applications. 

In this paper we single out the properties that are shared by various non-linear problems and study them  systematically.  
The property we investigate here is the linearity of the operators in the max-plus setting. Let us note that examples of  semigroups of max-plus linear  operators have already appeared in the literature (see e.g. \cite{AGL08, FK06, KM, Kr14, Mc06, Mc07}). However, there seems to be no systematic treatment of this class of evolution semigroups. The aim of this paper is therefore twofold. First, we give a proper definition of the max-plus linear semigroups of operators and present some  general properties.  Next, we show the usefulness of our object by pointing out some important  examples of such  semigroups.  We believe this can be a beginning of a fruitful theory. 

The paper is organized as follows. In Section 2 we define the max-plus vector spaces, max-additive and max-plus linear operators on them, as well as  max-additive and max-plus linear semigroups. In Section 3 we collect some general properties of such semigroups, such as contractivity in the Lipschitz norm. These properties are obtained by interlacing the max-plus and usual operations on function (max-plus) vector spaces and/or Banach lattices in an appropriate way.
In Section 4 we consider three important examples of nonlinear partial differential equations leading to max-additive and max-plus linear semigroups: the conservation law, Hamilton-Jacobi equation, and Hamilton-Jacobi-Bellman equation.  

\section{Definitions}
\subsection{The max-plus vector space}

First  let us denote by $\Rm$  the set $\mathbb{R} \cup \{-\infty\}$ equipped with the following two operations:
\begin{equation}\label{max-op}
a\oplus b=\max\{a,b\}\quad \text{ and } \quad a\otimes b = a+b.
\end{equation}
In particular we take 
$$a\oplus -\infty= a \;\;\mathrm{and}\;\; a\otimes -\infty =-\infty $$
for all $a \in \Rm$. 
The commutative idempotent semifield
$\Rm$ is usually called the \emph{max-plus algebra}.

In consistency with \cite{Mc06}, we call $(\cX,\oplus_{\cX},\otimes_{\cX})$  a \emph{max-plus vector space} over  $\Rm$ with zero element $0_{\cX}$ if for all $x,y,z\in\cX$ and $a,b\in \Rm $ it holds:
\begin{align*}
x\oplus_{\cX} y& = y \oplus_{\cX}x\in\cX, \quad a\otimes_{\cX} x\in\cX,\\
x\oplus_{\cX} 0_{\cX}  & =  x,\quad  a\otimes_{\cX} 0_{\cX} = 0_{\cX},\quad -\infty \otimes_{\cX} x=0_{\cX}, \quad 0\otimes_{\cX} x=x,\\
(x\oplus_{\cX} y)\oplus_{\cX} z &=x\oplus_{\cX} (y\oplus_{\cX} z ), \quad a\otimes_{\cX}(b\otimes_{\cX} x) = (a\otimes b)\otimes_{\cX} x,\\
 a\otimes_{\cX}(x\oplus_{\cX} y)& = (a\otimes_{\cX} x) \oplus_{\cX}(a\otimes_{\cX} y), \\
(a\oplus b)\otimes_{\cX} x &= (a\otimes_{\cX} x) \oplus_{\cX}(b\otimes_{\cX} x).
\end{align*}
Note that the last line implies $x\+x=x$ for all $x \in \cX$. {Max-plus vector spaces} are in the literature also known as moduloids, idempotent spaces or idempotent
semimodules (see e.g. ~\cite{LMS, KM, LM05}), where $\Rm$ can be replaced by a more general idempotent semifield or semiring. 

 A subset $\cY\subseteq\cX$ is called a \emph{max-plus subspace} of  max-plus vector space $\cX$, if $\cY$ is  invariant for the  operations $\oplus_{\cX} $ and $\otimes_{\cX} $. 

On every max-plus vector space $\cX$ we define the \emph{standard order} as the partial order $\preceq$ induced by $\oplus_{\cX}$:
\begin{equation}\label{order}
x\preceq y \iff x\oplus_{\cX} y=y.
\end{equation}

Typical examples of max-plus vector spaces are max-plus subspaces of a  max-plus vector space of functions
$$\cX=\{f\mid f\colon X \to \Rm\},$$
 where $X$ is a Haussdorff space or a $\sigma$-finite measure space, $\oplus_{\cX}$ a pointwise maximum, and $\otimes_{\cX}$ a pointwise addition, i.e.,
\begin{align*}
\left(f\oplus_{\cX} g\right)(x) &:= \max\left\{f(x),g(x)\right\}, \quad f,g\in\cX,\\
\left(a\otimes_{\cX} f\right)(x) &:= a+ f(x), \quad f\in\cX, \; a\in \Rm.
\end{align*}
In this case we define
\begin{equation}\label{zero}
0_{\cX}(x):= -\infty \quad\text{ and }\quad \Theta(x):=0, \quad x\in X.
\end{equation}
From now on we will consider $\cX$ to be such a function max-plus vector space and $\cC \subseteq \cX$ a max-plus subspace
or a subset invariant under $\oplus_{\cX}$.  We will denote 
$$\cCf=\{f \in \cC\mid f(x) > - \infty \;\;\mathrm{for} \;\;\mathrm{(almost)}\;\; \mathrm{all} \;\; x \in X \}.$$
Observe that if $\cC \subseteq \cX$ is a max-plus subspace, so is $\cCf \cup \{0_{\cX}\}$.

We will further assume that there exists a Banach space $\cE$ of functions  
from $X\to\mathbb{R}$ such that $\cCf \subseteq \cE$. The natural pointwise order in $\cE$ coincides on $\cCf$ with the standard order and therefore defines $\oplus_{\cX}$ on $\cCf$.  
For the sake of simplicity, we will omit the index $\cX$ in operations $\oplus$ and $\otimes$ whenever no confusion arises. 

We will interchangeably use the pointwise operations $+, -,  \cdot$ from  function space $\cE$ and  $\oplus, \otimes$ defined in $\cX$. Defining
\begin{equation*}
f^+:=f\oplus \Theta,\quad f^-:=-f\oplus \Theta \quad\text{and} \quad |f|:= f^+ + f^-
\end{equation*}
we obtain \emph{positive part}, \emph{negative part}, and \emph{absolute value} of every $f\in\cCf$ (the zero function $\Theta$ is defined in \eqref{zero}). Note that when applying these operations we ``move'' to the appropriate space where they are defined, hence $f^+$, $f^-$, $|f|$, or $\Theta$ are not necessarily contained in $\cCf$ or $\cX$.
We have the following relations 
\begin{align}
f &= f^+ - f^-,\label{f1} \\
f\oplus g &=  g + (f-g)^+ = f + (g-f)^+,\label{f2} \\
 |f-g|  &=  2 \cdot(f\oplus g) - f - g.\label{f3}
\end{align}

\subsection{Max-plus linear operators}
Let the subset $\cC\subseteq\cX$ be invariant under $\oplus$. Then an operator $T:\cC\to\cC$ is called
\begin{itemize}
\item[(i)] \emph{max additive}, if $T(f \oplus g)=Tf\oplus Tg$  for all  $f,g\in\cC$.
\end{itemize}
If, in addition, $\cC\subseteq\cX$ is also invariant under $\otimes$, then an operator $T:\cC\to\cC$ is called 
\begin{itemize}
\item[(ii)] \emph{plus homogeneous}, if $T(a\otimes f)=a\otimes Tf$  for all  $f\in\cC$  and $a\in \Rm$, and
\item[(iii)] \emph{max-plus linear}, if it is max additive and plus homogeneous.
\end{itemize}

For nonlinear operators the operator norm is not very convenient. Therefore we define the 
\emph{Lipschitz seminorm} of an operator $T:\cCf\to\cCf $ by
\begin{equation}\label{Lip-norm}
\|T\|_{\mathrm{Lip}} := \sup\left\{\frac{\|Tf-Tg\|}{\|f-g\|}\middle| f,g \in\cCf, f\neq g \right\}, 
\end{equation}
where $\|\cdot\|$ denotes the norm on $\cCf$ inherited from the Banach space $\cE$. Note that in case the operator $T$ is linear, $\|T\|_{\mathrm{Lip}}$ equals to the usual operator norm of $T$.

If the Lipschitz seminorms of the operators $T_1$ and $T_2$ exist, then it holds
\begin{equation}
\|T_1 T_2\|_{\mathrm{Lip}} \le \|T_1\|_{\mathrm{Lip}}\|T_2\|_{\mathrm{Lip}}\quad\text{and}\quad \|T_1 ^n\|_{\mathrm{Lip}}\le \|T_1\|_{\mathrm{Lip}}^n\quad \text{for }n\in\mathbb{N}.
\end{equation}

\subsection{Max-plus linear semigroup}

Let the subset $\cC\subseteq\cX$ be a max-plus vector subspace (or a subset invariant under $\oplus$).
We call a family $\Tt$ of max-plus linear (or max additive) operators on $\cC\subseteq\cX$  a one-parameter  \emph{max-plus linear (or max additive) semigroup} if $T(t)\cCf\subseteq \cCf$ and
\begin{equation}\label{FE}
\begin{array}{rcl}T(t+s) & = & T(t)T(s) \quad \text{ for all }t,s\ge0, \\
T(0) & = & Id_{\cC}.
\end{array}
\end{equation}

\begin{rem}\label{monotone} {\rm
It is easy to see that a max additive semigroup $\Tt$ is always \emph{monotone}, i.e.,  for every $f,g\in\cC$ the following implication holds:
\begin{equation*}
f \preceq g \Longrightarrow T(t)f\preceq T(t)g,
\end{equation*}
where $\preceq$ is the standard order on $\cX$, induced by the pointwise maximum.
}
\end{rem}

Let $\Tt$ be a semigroup of operators on $\cC$ such that  $T(t)\cCf \subseteq \cCf$ for all $t\ge 0$.
Then a semigroup $\Tt$   is called \emph{strongly continuous} on $\cC$, if the orbit mappings 
\begin{equation}\label{SC}
t\mapsto T(t)f : \RR_+ \to \cCf
\end{equation}
are continuous  for every $f\in\cCf$.

The  \emph{infinitesimal generator} $A$ of a strongly continuous  semigroup $\Tt$ is defined as 
\begin{equation}\label{generator}
Af:=\lim_{t\downarrow 0}\frac{T(t)f-f}{t} ,
\end{equation}
where the domain $D(A)$ is the set of all $f\in\cCf$ for which the above limit exists. 
Note that in the theory of nonlinear semigroups one can in general consider also multivalued operators (see \cite{CL71,GN14}). However,  because of the definition \eqref{generator}, we treat only single valued operators.

\section{Properties of max-plus linear semigroups}

\subsection{Standard constructions} 

There are several ways to construct new strongly continuous semigroups from a given one. 
Let us check which constructions preserve max-plus linearity. 

In all following lemmata we assume that $\Tt$ is a max-plus linear strongly continuous semigroup on a max-plus subspace $\cC\subseteq\cX$.
The proofs are all straightforward and therefore omitted (see also \cite[Sec. I.5.b]{EN00}). 

\begin{lem}[Similarity]\label{similar}
Let $\cY$ be a max-plus vector space and $V:\cY\to\cX$ a max-plus linear isomorphism such that 
$ V^{-1} (\cCf)=V^{-1} (\cC )_{\mathrm{fin}} $. Then
\begin{equation*}
S(t):= V^{-1} T(t) V,\quad t\ge0,
\end{equation*}
form a max-plus linear strongly continuous semigroup on $V^{-1}(\cC)\subseteq\cY$. \end{lem}

\begin{lem}[Restriction]\label{subspace}
Let $\cD\subseteq\cC$ be a max-plus subspace such that $\cDf$ is closed in $\cE$,  $T(t)\cD\subseteq \cD$ and 
$T(t)\cDf\subseteq \cDf$  for all $t\ge 0$. Then 
the restrictions
\begin{equation*}
S(t):= T(t)|_{\cD} ,\quad t\ge0,
\end{equation*}
form a max-plus linear strongly continuous semigroup on $\cD$.
\end{lem}

For a max-linear subspace $\cD\subset\cC$  we define an equivalence relation on $\cC$ (see e.g. \cite{KP, KM})
 by
\begin{equation*}
f_1\sim f_2 \iff \text{ there are } g_1,g_2\in\cD \text{ such that } f_1 \oplus g_1 = f_2 \oplus g_2.
\end{equation*}
We denote the appropriate equivalence classes by $[f]_{\cD}$, $f\in\cC$, and the set of this classes by
\begin{equation*}
\cC/\cD: = \{[f]_{\cD}\mid f\in \cC\}.
\end{equation*}
The operations $\oplus$ and $\otimes$ on $\cC/\cD$ are defined naturally by
\begin{align*}
[f]_{\cD} \oplus [g]_{\cD} &:= [f\oplus g]_{\cD}, \quad f,g\in\cC,\\
a\otimes [f]_{\cD} &:= [a\otimes f]_{\cD},\quad f\in\cD, a\in\Rm.
\end{align*}
It is not difficult to  verify that  these two operations  make $\cC/\cD$  a max-linear vector space over $\Rm$ (with zero element $\cD$), which we call the \emph{max-plus quotient space}. 

\begin{lem}[Quotient]\label{quotient}
Let $\cD\subset\cC$ be a max-linear subspace that is invariant for $\Tt$ and $\cC/\cD$  a max-plus quotient space.
Then 
\begin{equation*}
S(t) [f]_{\cD} := \left[T(t)f\right]_{\cD},\quad f\in\cC, t\ge0,
\end{equation*}
defines a max-plus linear semigroup on $\cC/\cD$.
\end{lem}

\begin{lem}[Product]\label{product}
Let $\left(U(t)\right)_{t\ge 0}$ be a max-plus linear  strongly continuous semigroup on $\cC$  such that $T(t)U(t) = U(t)T(t)$ for all $t\ge 0$. Then the products
\begin{equation*}
S(t):= T(t)U(t),\quad t\ge0,
\end{equation*}
form a max-plus linear strongly continuous semigroup on $\cC$.
\end{lem}

\begin{lem}[Rescaling]\label{rescaling}
Let $\alpha >0, \beta\in\mathbb{R}$. Then 
\begin{equation*}
S(t):= e^{\beta t} T(\alpha t),\quad t\ge0,
\end{equation*}
form  max-plus linear strongly continuous semigroup on $\cC$, the so-called \emph{rescaled semigroup}.
\end{lem}

\subsection{Contraction property} 

Let $\Tt$ be a strongly continuous semigroup of operators on a subset $\cC\subseteq\cX$ such that $\cCf$ is closed 
in $\cE$. In the literature treating the nonlinear operator semigroups  (see e.g. \cite{CL71,GN14}) it is usually assumed in advance that   for some  $\omega\in\RR$,
\begin{equation*}
\|T(t)f-T(t)g\|\le e^{\omega t}\|f-g\|\quad \text{for all }t\ge 0 \text{ and }f,g\in\cC,
 \end{equation*}
that is,
\begin{equation}
 \|T(t)\|_{\mathrm{Lip}}\le e^{\omega t}\quad \text{for all } t\ge 0.
 \end{equation}
If $\omega=0$, $\Tt$ is called a \emph{contraction semigroup}.
We will see that this condition is always satisfied for a max-plus linear semigroup on certain function spaces.

Let $\Omega$ be some measure space with a $\sigma$-finite measure $\mu$. The spaces $\llp(\Omega,\mu) \cup \{-\infty\}$, $1\le p<\infty$, are invariant for the max operation $\oplus$, and 
$\lli(\Omega,\mu) \cup \{-\infty\}$ is a  max-plus vector space. Here $-\infty$ denotes the equivalence class of functions almost everywhere equal to $-\infty$.  Inspired by \cite{CT80} we observe the following. 

\begin{prop}\label{cons-L1}
Let $\Tt$ be a max additive semigroup on $\cC\subseteq \ll1(\Omega,\mu)\cup \{-\infty\}$, where the set $\cC$ is invariant for the max operation.  Assume moreover, that 
\begin{equation*}
\int_{\Omega} T(t)f d\mu = \int_{\Omega}f d\mu\quad \text{for all } t\ge 0\text{ and } f\in\cCf.
\end{equation*}
Then every $T(t)$ is an isometry and thus
\begin{equation*}
 \|T(t)\|_{\mathrm{Lip}} = 1\quad \text{for all } t\ge 0.
 \end{equation*}
\end{prop}
\begin{proof}
Using relation \eqref{f3} and max additivity of the semigroup we obtain
\begin{equation*}
\left|T(t)f-T(t)g\right| = 2\cdot  T(t)(f\oplus g) - T(t)g - T(t)f
\end{equation*}
for $f,g \in \cCf$.
Since by assumption $T(t)$ preserves the integral it follows by  \eqref{f3}
\begin{align*}
\|T(t)f - T(t)g\|_1 &= \int_{\Omega} \left| T(t)f - T(t)g\right| d\mu \\
&=  2\int_{\Omega} T(t)(f\oplus g) d\mu - \int_{\Omega} T(t)g d\mu-\int_{\Omega}  T(t)f d\mu\\
&=  \int_{\Omega} \left(2 (f\oplus g)- g - f\right) d\mu\\
&= \int_{\Omega} \left| f - g\right| d\mu = \|f - g\|_1
\end{align*}
for all $f,g\in\cCf$.
\end{proof}

In case of $\lli(\Omega,\mu)$ max-plus linearity of the semigroup yields similar conclusion.

\begin{prop}\label{cons-Linfty}
Let $\Tt$ be a max-plus linear semigroup on a max-plus subspace $\cC\subseteq \lli(\Omega,\mu) \cup \{-\infty \}$.  Then 
\begin{equation*}
 \|T(t)\|_{\mathrm{Lip}}\le1\quad \text{for all } t\ge 0.
 \end{equation*}
\end{prop}
\begin{proof}
First observe that by (\ref{f2})
\begin{equation*}
f\oplus g \preceq g+\|(f-g)^+\|_{\infty} \quad\text{as well as}\quad f\oplus g \preceq f+\|(g-f)^+\|_{\infty}
\end{equation*}
for $f,g \in \cCf$.
Plugging this into \eqref{f3} and using monotonicity (see Remark \ref{monotone}) and  plus homogeneity of the semigroup yields
\begin{align*}
\left|T(t)f\right. - \left. T(t)g\right| &= 2\cdot T(t)\left(f\oplus g\right) - T(t)g- T(t)f\\
 &\preceq  
 T(t)\left(g+\|(f-g)^+\|_{\infty}\right) - T(t)g \\
 & \quad + T(t)\left(f+\|(g-f)^+\|_{\infty}\right) - T(t)f \\
&= \|(f-g)^+ \|_{\infty} + \|(g-f)^+\|_{\infty} = \|f-g\|_{\infty}
\end{align*}
for all $f,g\in\cCf$.
\end{proof}

Let  $\cb(X)$ be a Banach space of all  bounded continuous functions (on a Haussdorff space $X$) equipped with the supremum norm. A similar proof as above proves the following result.

\begin{prop}\label{cons-ck}
Let $\Tt$ be a max-plus linear semigroup on a max plus subspace $\cC\subseteq  \cb(X) \cup \{-\infty\}$. Then 
\begin{equation*}
 \|T(t)\|_{\mathrm{Lip}}\le1\quad \text{for all } t\ge 0.
 \end{equation*}
\end{prop}

\subsection{The generator}

The infinitesimal generator of a  strongly continuous semigroup of \emph{linear} operators is always linear. The generator of a max-plus linear semigroup is however never max-plus linear. Observe that from the definition \eqref{generator} of the generator  it follows 
that it is always translation invariant, i.e.,  
$$A(a\otimes f ) = Af\quad\text{for all }f\in\cCf, a\in\Rm,$$
so $A$ is never plus homogeneous. It is also not always max additive, as the following example shows.

\begin{ex}
{ \rm Let $\cC=\buc(\mathbb{R}) \cup \{-\infty\}$ be the space of bounded uniformly continuous functions on $\mathbb{R}$ with the function $-\infty$ appended. Note that it is a function max-plus vector space for the pointwise operations. Take the left translation semigroup  on $\cC$ defined as 
$$T(t)f(s):=f(s+t),\quad f\in\cC.$$ It is well-known that $\Tt$ is a strongly continuous semigroup of linear operators (cf.~\cite[Sec.~I.4.c]{EN00}). The operators $T(t)\colon\cC\to\Rm$ are also max-linear:
\begin{align*}
\left(T(t)f\oplus T(t)g)\right)(s)&=\max\{T(t)f(s),T(t)g(s)\}=\max\{f(s+t),g(s+t)\}\\
&=(f\oplus g)(s+t) = \left(T(t)(f\oplus g)\right)(s)
\end{align*}
and
\begin{align*}
T(t)(a\otimes f)(s) 
&= (a\otimes f)(s+t) = a+ f(s+t) \\
&= a + T(t)f(s) = \left(a\otimes T(t)f\right)(s). 
\end{align*}
The infinitesimal generator of $\Tt$ is the operator (see \cite[Sec.~II.2.b]{EN00})
$$Af:=f',\quad D(A)=\buc(\mathbb{R})\cap \c1(\mathbb{R}).$$ 
Operator $A$ is linear, but not max-additive. E.g., taking $f(x)=e^{-2x^2}$ and $g(x)=e^{-x^2}$ we have $f\oplus g = g$ and thus $A(f\oplus g)=Ag$ while it is easy to see that $Af\oplus  Ag \ne Ag$.
}
\end{ex}

In generation theorems for nonlinear semigroups (e.g., Crandall-Ligget Theorem \cite{CL71}) a necessary condition for the generator $A$ is \emph{dissipativity}, i.e., for each $\alpha>0$ the inverse operator $(I-\alpha A)^{-1}$ exists and $\|(I-\alpha A)^{-1}\|_{\mathrm{Lip}}\le 1$. 
By \cite[Theorem III.1.1]{B76}, Proposition \ref{cons-Linfty} and Proposition \ref{cons-ck} the following holds in our case.

\begin{cor}\label{diss} Let $\cE$ be  $\lli(\Omega,\mu)$ or  $\cb(X)$.  Let $\Tt$ be a strongly continuous max-plus linear semigroup on a  max-plus linear subspace $\cC\subseteq \cE \cup \{-\infty\}$ such that $\cCf$ is closed in $\cE$.
Then its infinitesimal generator is dissipative.
\end{cor}

\section{Examples of max additive and max-plus linear semigroups\protect\footnote{After this paper has already been published online it has been noticed by B. Adreianov that the
conditions of Propositions \ref{CL-semigroup} and \ref{HJ-semigroup} are not stated accurately enough and the statements do
not hold in the present form. This has been corrected in the enclosed Erratum.}}

Here we consider some important examples of nonlinear semigroups and show that they are max-additive or max-plus linear.

\subsection{Scalar Conservation Law}

We are interested in the solutions to the following Cauchy problem 
\begin{equation}\label{eq:cl}\tag{CL}
\left\lbrace\begin{aligned}
u_t + f(u)_x &= 0,\quad t>0, \;x\in\mathbb{R}, \\
u(x,0)&=h(x),\quad x \in\mathbb{R}.
\end{aligned}\right. 
\end{equation}
The quasilinear equation in the first line of \eqref{eq:cl} is known as the (scalar) \emph{conservation law}.
Even for smooth initial conditions, the classical (continuously differentiable) solutions to this problem do not always exist therefore we need to generalize the concept of solution. 
A $\ll1$-function $u$ on $\mathbb{R}\times[0,\infty)$ is called a \emph{weak solution} to  \eqref{eq:cl} if 
\begin{equation}\label{weak}
\int_{0}^{\infty} dt \int_{-\infty}^{\infty}  \left(u \cdot \psi_t + f(u)\cdot \psi_x\right)\; dx +  \int_{-\infty}^{\infty} h(x)\cdot \psi(x,0)\; dx =0
\end{equation}
holds for every $C^1$-function $\psi$ on $\mathbb{R}\times[0,\infty)$ with compact support.
Weak solutions are not unique and, in order to obtain the physically correct solution, one has to impose the right entropy condition. 
There is a rich mathematical theory on this topic, see for example \cite{Br05,Cr72} or \cite[Sec.4]{GN14}.

We will use the compact description of the right solutions due to Kru\v{z}kov. Assume that $f:\mathbb{R}\to\mathbb{R}$ is locally Lipschitz continuous and that $h\in \ll1(\mathbb{R})$. 
The \emph{entropy solution} to \eqref{eq:cl} is a continuous map $v\colon [0,\infty)\to\ll1(\mathbb{R})$ which satisfies $v(0)=h$ together with
\begin{equation}\label{kruzkov}
\int_{0}^{\infty} dt \int_{-\infty}^{\infty}  \big( |v-k | \cdot\psi_t + (f(v)-f(k))\: \sgn(v-k)\cdot \psi_x\big)\; dx \ge 0
\end{equation}
for every $k\in \mathbb{R}$ and every non-negative function $\psi\in\cc1(\mathbb{R}^2)$, whose compact support is contained in the half plane where $t>0$.
Then the unique bounded entropy solution to \eqref{eq:cl} are given as trajectories $t\mapsto T(t)h$ where $\Tt$ is a strongly continuous semigroup (see \cite[Theorem 6.3]{Br05}). 
The operators in this semigroup are known to be nonlinear. However, it is not difficult to see that they are max-additive.

\begin{prop}\label{CL-semigroup}
The semigroup $\cT^{\mathrm{CL}}:=\Tt$,  where  $u(t,x)=T(t)h(x)$ is the unique entropy 
solution to \eqref{eq:cl} and where $T(t)(-\infty):= -\infty$, is a max additive strongly continuous semigroup on $\ll1(\mathbb{R}) \cup \{-\infty\}$. 
\end{prop}

\begin{proof}
Our aim is to show that 
$$T(t)h_1\oplus T(t) h_2 = T(t) (h_1\oplus h_2) 
\text{ for any }h_1,h_2\in\ll1(\mathbb{R}).$$ 
Let $\psi\in\cc1(\mathbb{R}^2)$ be any appropriate test function and denote by
$\psi^1$  the smooth cutoff function that coincides with $\psi$ on the compact set 
$${ \{ (t,x)\mid T(t) h_1(x)\geq T(t) h_2(x) \}}\cap\supp \psi$$ and equals $0$ outside of a neighbourhood of this set.
Analogously define  the smooth cutoff function $\psi^2$ that coincides with $\psi$ on the complementary set and so $\supp \psi^1 \cup \supp \psi^2 = \supp \psi$.
 Then we can write
{\small
\begin{align*}
&\int_{0}^{\infty} dt \int_{-\infty}^{\infty}  \bigg( \big|T(t) h_1\oplus T(t) h_2 - k\big| \cdot\psi_t \\
& +  \Big(f\big(T(t) h_1\oplus T(t) h_2\big)-f(k)\Big)\: \sgn\big( T(t) h_1\oplus T(t) h_2 - k\big) \cdot\psi_x\bigg)\; dx
  \\
&=\int_{0}^{\infty} dt \int_{-\infty}^{\infty}  \bigg( |T(t) h_1 - k | \cdot\psi^1_t + 
 \left(f(T(t) h_1)-f(k))\: \sgn( T(t) h_1 - k\right) \cdot\psi^1_x\bigg)\; dx
  \\
&+  \int_{0}^{\infty} dt \int_{-\infty}^{\infty}  \bigg( |T(t) h_2 - k | \cdot\psi^2_t + 
 \left(f(T(t) h_2)-f(k))\: \sgn( T(t) h_2 - k\right) \cdot\psi^2_x\bigg)\; dx\\
 \end{align*}
 }
Since $T(t)h_1$ and $T(t)h_2$ are the entropy solutions to \eqref{eq:cl} with initial conditions $h=h_1$ and $h=h_2$, respectively, both integrals above are greater or equal than $0$ for every $k\in\mathbb{R}$. Moreover, 
$T(t) h_1\oplus T(t) h_2\vert_{t=0} = h_1\oplus  h_2 $, hence $T(t) h_1\oplus T(t) h_2$ is the entropy solution to \eqref{eq:cl} with initial condition $h=h_1\oplus h_2$. By uniqueness of the entropy solutions we obtain 
$T(t)h_1\oplus T(t) h_2 = T(t) (h_1\oplus h_2)$.
\end{proof}

Max-additivity directly leads to monotonicity of the semigroup $\cT^{\mathrm{CL}}$ defined in Proposition \ref{CL-semigroup}. Furthermore, the solutions $u(t,x)=T(t)h(x)$ to \eqref{eq:cl} preserve the integral, i.e.,
\begin{equation*}
\int_{-\infty}^{\infty}  T(t)h(x) \; dx =\int_{-\infty}^{\infty} h(x)\; dx,\quad t\ge 0.
\end{equation*}
Hence, by Proposition \ref{cons-L1}, $\cT^{\mathrm{CL}}$ consists of isometries and  $\|\cT^{CL}\|_{\mathrm{Lip}} = 1$ (compare with \cite[Theorem 6.3.(ii)]{Br05}). Moreover, by  \cite[Theorem III 1.1]{B76} the infinitesimal generator of $\cT^{CL}$ is dissipative.

Note that the semigroup $\cT^{\mathrm{CL}}$ is in general not plus homogeneous. This happens only if the function $f$ appearing in (CL) is nice enough.

\subsection{Hamilton-Jacobi Equation}

Let us now consider a similar  initial value problem
\begin{equation}\label{eq:HJ}\tag{HJ}
\left\lbrace\begin{aligned}
u_t + f(\nabla u) &= 0,\quad t>0, \;x\in\mathbb{R}^n, \\
u(x,0)&=h(x),\quad x \in\mathbb{R}^n.
\end{aligned}\right. 
\end{equation}
Here $\nabla u$ denotes the gradient with respect to $x\in\mathbb{R}^n$ and $f,h\colon\mathbb{R}^n\to \mathbb{R}$ are given functions.
Writing $f(\nabla u) = H(x,\nabla u)$, we call $H\colon\mathbb{R}^n\times \mathbb{R}^n\to \mathbb{R}$ the \emph{Hamiltonian function} and  \eqref{eq:HJ}  the \emph{Hamilton-Jacobi equation}, which often appears in optimization problems.

As before, the classical solution to \eqref{eq:HJ} does not always exist and one again has to consider some generalized solutions that correspond to previously defined Kru\v{z}kov entropy solutions. They are called \emph{viscosity solutions} and are defined in \cite{CEL84} where also the existence and uniqueness result is proved assuming that $H\in C\left(\mathbb{R}^n\right)$ and $h\in \buc\left(\mathbb{R}^n\right)$.

\begin{prop}\label{HJ-semigroup}
The semigroup $\cT^{\mathrm{HJ}}:=\Tt$,  where  $u(t,x)=T(t)h(x)$ is the unique viscosity
solution to \eqref{eq:HJ} and where  $T(t)(-\infty):= -\infty$, is a strongly continuous max-plus linear semigroup on $\buc\left(\mathbb{R}^n\right) \cup \{-\infty\}$. 
\end{prop}
\begin{proof}
For max additivity of the semigroup we refer to \cite[Prop.~1.3.(a)]{CEL84}.

Now let $u(t,x)=T(t)h(x)$ be the unique viscosity solution to \eqref{eq:HJ}, $a\in \Rm$, and 
\begin{equation*}
\widetilde{u}(t,x) := \left(a\otimes T(t)h\right)(x) = a + T(t)h(x)  = a+ u(t,x). 
\end{equation*}
Then $\widetilde u_t = u_t$ and  $\nabla \widetilde u= \nabla u$,
 therefore $\widetilde u$ solves \eqref{eq:HJ} with the initial function $\widetilde u (0,x)= \left(a\otimes h\right)(x)$. By the uniqueness of the viscosity solutions we have
\begin{equation*}
\widetilde{u}(t,x) = \left(a\otimes T(t)h\right)(x)=T(t)\left(a\otimes h\right)(x) \text{ for all }x\in\mathbb{R}^n,
\end{equation*}
hence $\Tt$ is also plus homogeneous.
\end{proof}

By Proposition \ref{cons-ck}, the max-plus linear semigroup $\cT^{\mathrm{HJ}}$ consists of contractive operators and by Corollary \ref{diss} its generator is dissipative (compare with \cite[Remark 2.1, Proposition 5.1]{CEL84}).

\subsection{Hamilton-Jacobi-Bellman Semigroup}

Consider the following deterministic finite time-horizon optimal control problem (see \cite{AGL08,FS93}, or in a more general setting \cite{LN82}):
\begin{subequations}\label{eq:FHOCP}
\begin{align}
\text{maximize} & \int_0^T\ell(x(s),u(s))\,ds+\phi(x(T))\\
\text{ such that}&\notag\\
\dot{x}(s)=&f(x(s),u(s)),\quad x(0)=x
\label{eq:1}
\end{align}
\end{subequations}
with $x(s)\in X,$ $u(s)\in U,$ $0\le s\le T.$
Here the \emph{state space} $X$ is a subset of $\RR^n$, the set of \emph{control values} $U$ is a subset of $\RR^m$, the time-horizon is $T>0$. The \emph{initial condition} $x\in X$ is given, the control function $u(\cdot)$ is bounded and Lebesgue-measurable, the map $x(\cdot)$ is absolutely continuous. We assume that the \emph{instantaneous reward} $\ell:X\times U\to\RR$ and the \emph{dynamics} $f:X\times U\to\RR^n$ are sufficiently regular maps, and the \emph{terminal reward} is a map $\phi:X\to\RR\cup\{-\infty\}$.

The so-called \emph{value function} $v$ associates to any $(x,t)\in X\times [0,T]$ the following supremum
\begin{equation}\label{eq:value}
v(x,t):=\sup \int_0^t\ell(x(s),u(s))\,ds+\phi(x(t)),
\end{equation}
where the supremum is taken under the constraint \eqref{eq:1}.\\
According to \cite[Theorem I.5.1]{FS93}, the value function $v$ is the solution of the \emph{Hamilton-Jacobi-Bellman partial differential equation}:
\begin{equation}\label{eq:HJB}\tag{HJB}
 \left\{\begin{aligned}
-v_t+H(x,\nabla v)&=0,\quad (x,t)\in X\times (0,T];\\
v(x,0)&=\phi(x),\quad x\in X,
\end{aligned}\right.
 \end{equation}
where
\begin{equation}\label{eq:Ham}
H(x,p)=\sup_{u\in U}\left(\ell(x,u)+p\cdot f(x,u)\right),\quad (x,p)\in X\times\RR^n
\end{equation}
for given $U\subset\RR^m$ is the \emph{Hamiltonian} of the problem. This system is clearly a special case of the Hamilton-Jacobi equation \eqref{eq:HJ}.

For a fixed $\phi\in \buc(X)$, denote the maps
\begin{equation}\label{eq:semigroup}
S(T)\phi:=v(\cdot,T),
\end{equation}
that associate to any $\phi$  -- the terminal reward of \eqref{eq:FHOCP} -- the value function \eqref{eq:value} on horizon $T.$ Since $\phi$ is the initial condition and $v$ the solution of \eqref{eq:HJB}, $\cT^{\mathrm{HJB}}:=(S(t))_{t\geq 0}$ is the \emph{evolution semigroup of} \eqref{eq:HJB} (the so-called \emph{Lax-Oleinik semigroup}).

Due to \cite[Example II.3.1]{FS93} and \cite{LN82}, the family of maps $\cT^{\mathrm{HJB}}$ forms a one-parameter \emph{strongly continuous semigroup on $\buc(X)$}, consisting of nonlinear and monotone
operators. By Proposition \ref{HJ-semigroup} it is also a max-plus linear semigroup on $\buc(X)\cup\{-\infty\}$. This property for the semigroup was already observed by Maslov \cite{Mas73}.

As semigroup $\cT^{\mathrm{HJ}}$, also  $\cT^{\mathrm{HJB}}$ is a contraction semigroup whose generator is dissipative (see Proposition \ref{cons-ck} and Corollary \ref{diss}).

\begin{rem}{\rm 
To the contrary of the general Hamilton-Jacobi equation \eqref{eq:HJ}, here the semigroup $\cT^{HJB}=\left(S(t)\right)_{t\geq 0}$ has an ``explicit form''
\begin{equation}\label{eq:semigroupexpl}\left(S(t)\phi\right)(x)=v(x,t)=\sup \int_0^t\ell(x(s),u(s))\,ds+\phi(x(t)),\end{equation}
where the supremum is taken under the constraint \eqref{eq:1}. If $f$ is locally Lipschitz in $x$ (i.e., in the first variable), then \eqref{eq:1} has a unique solution $x(\cdot)$ for each initial data $x\in X$ and $u(\cdot)\in \lli\left([0,T],U\right)$. Hence, the supremum in \eqref{eq:value} or in \eqref{eq:semigroupexpl} can be actually taken (for a fixed $x\in X$) in $u(\cdot)\in \lli\left([0,T],U\right)$.
}
\end{rem}

\begin{rem}{\rm The generator $(A,D(A))$ of the semigroup $\cT^{\mathrm{HJB}}$ is actually the Hamilton-Jacobi-Bellman operator $H(x,\nabla \phi)$, that is
\begin{align*}
(A\phi)(x)&=\lim_{t\to 0}\frac{1}{t}(S(t)\phi-\phi)(x)=H(x,\nabla\phi)=\sup_{u\in U}\left(\ell(x,u)+\nabla\phi\cdot f(x,u)\right)\\
D(A)&=\left\{\phi\in \buc(X):\nabla\phi\in \buc(X)\right\},
\end{align*}
see \cite[II. (3.13)]{FS93} or \cite{LN82}.}
\end{rem}

{\it Acknowledgments.}
The authors thank  J.A. Goldstein, R. Nagel and M. Kandi\'{c} for useful comments. 
The first and the second author were supported in part by the Slovenian Research Agency. \\

\medskip
{\small

\noindent
\emph{Marjeta Kramar Fijav\v{z}}, University of Ljubljana, Faculty of Civil and Geodetic Engineering, Jamova 2, SI-1000 Ljubljana, Slovenia /
Institute of Mathematics, Physics, and Mechanics,
Jadranska 19, SI-1000 Ljubljana, Slovenia,\\
\texttt{marjeta.kramar@fgg.uni-lj.si}

\medskip
\noindent
\emph{Aljo\v{s}a Peperko},
University of Ljubljana, Faculty of Mechanical Engineering, A\-\v{s}ker\-\v{c}e\-va 6, SI-1000 Ljubljana, Slovenia /
Institute of Mathematics, Physics, and Mechanics,
Jadranska 19, SI-1000 Ljubljana, Slovenia, \\\texttt{aljosa.peperko@fs.uni-lj.si}

\medskip
\noindent
\emph{Eszter Sikolya},
ELTE TTK, Department of Applied Analysis and Computational Mathematics, 
P\'azm\'any P\'eter s\'et\'any 1/C, H-1117 Budapest, Hungary,\\
\texttt{seszter@cs.elte.hu}}
\end{document}


\maketitle

The examples of max-additive and max-plus linear semigroups in the last section of the paper \cite{KPS15} are given inaccurately.
\cite[Proposition 4.1]{KPS15} is not true as stated and \cite[Proposition 4.2]{KPS15} does not hold without some additional assumptions.

In the case of conservation law $(CL)$ studied in \cite[Subsec.~4.1]{KPS15}, a smaller set of weak solution should be considered.  
A function  $v\colon[0,\infty)\to\ll1(\mathbb{R})$ is called an  \emph{isentropic solution} to $(CL)$, if it satisfies the initial condition $v(0)=h$ and  the 
Kru\v{z}kov condition $(14)$ from \cite{KPS15} with the equality sign. 

We now give the proper formulation of \cite[Proposition 4.1]{KPS15} and its proof.
Let 
\begin{equation*}
\begin{split}\cC:=&\{h\in\ll1(\mathbb{R})\mid (CL) \text{ admits an isentropic solution } v \text{ with }v(0)=h\}\\
&\cup \{-\infty\}.
\end{split}
\end{equation*}
\begin{pr}[{\cite[Proposition 4.1]{KPS15}}]
The semigroup $\cT^{\mathrm{CL}}:=\Tt$,  where  $u(t,x)=T(t)h(x)$ is the unique isentropic
solution to $(CL)$ and where \\$T(t)(-\infty):= -\infty$, is a max additive strongly continuous semigroup on $\cC$. 
\end{pr}

\begin{proof}
By \cite[Theorem 1]{KP07}, $T(t)h_1\oplus T(t)h_2$ is an isentropic solution to (CL) for any $h_1,h_2\in\cC$.
Moreover, 
$T(t) h_1\oplus T(t) h_2\vert_{t=0} = h_1\oplus  h_2 $, hence $T(t) h_1\oplus T(t) h_2$ is the isentropic solution to (CL) with initial condition $h=h_1\oplus h_2\in\cC$. By uniqueness of the isentropic solutions we obtain 
$T(t)h_1\oplus T(t) h_2 = T(t) (h_1\oplus h_2)$.
\end{proof}

\begin{rmk}
The isentropic solutions to $(CL)$ often exist only on some limited time interval. However, our results also hold in this case considering the local time flows  $\left(T(t)\right)_{t\in[0,T]}$ defined on
\begin{equation*}
\begin{split}
\cC_T:=\{&h\in\ll1(\mathbb{R})\mid (CL) \text{ admits an isentropic solution } v \text{  on } [0,T] \\
&\text{ with }v(0)=h\} \cup \{-\infty\}.
\end{split}
\end{equation*}
\end{rmk}

In the case of Hamilton-Jacobi equation $(HJ)$ considered in \cite[Subsec.~4.2]{KPS15} we should specify the choice of the generalized  solutions and regularity assumptions.  We assume the following for the function $f$ appearing in the problem $(HJ)$.
\begin{enumerate}
\item For every $x\in\mathbb{R}^n$, $\|x\|=1$, there exist the limit $\lim_{r\downarrow 0} r f(x/r)$.
\item For any $(x,r), (x',r')\in\mathbb{R}^n\times\mathbb{R}_+$ with $\|x\|^2 +r^2 \le 1$ and $\|x'\|^2 +r'^2 \le 1$,
$$\left| r f(x/r)-r' f(x'/r')\right|\le K \left(\|x-x'\|^2 + (r-r')^2\right)^{1/2},$$
where $K>0$ is some constant.
\item The function $f$ is convex.
\end{enumerate}
In the literature there are many notions for generalized weak solutions to $(HJ)$: Crandall-Lions' viscosity solutions, minimax solutions by Subbotin, Maslov idempotent weak solutions, Kru\v{z}kov generalized solutions.  However, under our assumptions they all agree, see \cite[Sec. 5]{Ro05}, and 
\cite[Proposition 4.2]{KPS15} holds as stated in the paper. 
For the proof we refer directly to \cite[Theorem 3.2]{Ro05} (note that there the min-plus terminology is used, therefore the concavity instead of convexity of $f$ is assumed). 

\begin{rmk}
Stability of solutions of $(HJ)$ (understood in the viscosity sense of Crandall-Lions)
under the $\oplus$-operation is closely related to the properties of $\mathrm{liminf}$
of a sequence of viscosity solutions of first order Hamilton-Jacobi equations
with convex Hamiltonians, see \cite{BJ90} and \cite[Th. 2.1 and 2.3]{Ba93}.
The fact that the maximum of two viscosity solutions is a viscosity sub-solution
is classical, see \cite{CIL92}.
Heuristically, the fact that the maximum of two viscosity solutions
is also a viscosity super-solution stems from the semi-concavity property
of viscosity solutions, see \cite[p.1125]{Ba93}.
\end{rmk}
\begin{rmk}
We add a reference for yet another proof of the max-linearity of  Hamilton-Jacobi-Bellman semigroup presented in Subsection 4.3. In \cite{DD99} this fact is proved using a probabilistic approach. 
\end{rmk}

\noindent{\it Acknowledgments.}
The first author thanks Guy Barles for fruitful discusions.

\medskip
{\small

\noindent
\emph{Boris Andreianov},
LMPT CNRS UMR7350, Universit\'e de Tours,
Parc Grandmont, 37200 Tours, France,
\texttt{boris.andreianov@lmpt.univ-tours.fr}

\medskip
\noindent
\emph{Marjeta Kramar Fijav\v{z}}, University of Ljubljana, Faculty of Civil and Geodetic Engineering, Jamova 2, SI-1000 Ljubljana, Slovenia /
Institute of Mathematics, Physics, and Mechanics,
Jadranska 19, SI-1000 Ljubljana, Slovenia,\\
\texttt{marjeta.kramar@fgg.uni-lj.si}

\medskip
\noindent
\emph{Aljo\v{s}a Peperko},
University of Ljubljana, Faculty of Mechanical Engineering, A\-\v{s}ker\-\v{c}e\-va 6, SI-1000 Ljubljana, Slovenia /
Institute of Mathematics, Physics, and Mechanics,
Jadranska 19, SI-1000 Ljubljana, Slovenia, \\\texttt{aljosa.peperko@fs.uni-lj.si}

\medskip
\noindent
\emph{Eszter Sikolya},
ELTE TTK, Department of Applied Analysis and Computational Mathematics, 
P\'azm\'any P\'eter s\'et\'any 1/C, H-1117 Budapest, Hungary,\\
\texttt{seszter@cs.elte.hu}}